\documentclass[12pt,a4paper]{article}
\usepackage{amsmath}
\usepackage{amsfonts}
\usepackage{amssymb}
\usepackage{amsmath,authblk,hyperref}
\newcommand{\RNum}[1]{\uppercase\expandafter{\romannumeral #1\relax}}

\title{\textbf{ Common Fixed Point Theorems for Four Transformations in Vector $S$-metric Spaces }}
\author[1]{Pooja Yadav\thanks{poojayadav.math.rs@igu.ac.in}}
\author[2]{Mamta Kamra\thanks{mkhaneja15@gmail.com}}
\affil[1, 2]{Department of Mathematics, Indira Gandhi University Meerpur(Rewari), Haryana-122502, India.}
\date{}
\begin{document}
\maketitle
\begin{abstract}
\noindent We establish some common fixed point results for four transformations in vector $S$-metric spaces by using the notion of weakly compatibility (WC) and occasionally weakly compatibility (OWC). The first theorem is proved by using the concept of CLR$_{(p, q)}$ property (common limit range property w.r.t. transformations $p$ and $q$) and weakly compatiblity whereas second theorem is proved by using the concept of occasionally weakly compatiblity.
\end{abstract}
\section{Introduction}
The fixed point theory is used in control theory, optimization theory, and many branches of analysis. Initially, Sessa \cite{Sessa} proved common fixed point results for weakly commutative pair of transformations in 1992. In 1996, Jungck \cite{Jungck} gave the notion of weakly compatible transformations. In 2002, Aamir and Moutawakil \cite{Aamri} explained the concept of (E. A.) property which needs the closedness (or completeness) of the subspaces for the existence of a common fixed point. Sintunavarat et al. \cite{Sintunavarat} gave the concept of $CLR_Q$ in the metric space for the pair of transformations in 2011 which relaxes the closedness (or completeness) of the subspaces. Many researchers [\cite{Roldan}, \cite{Manro}, \cite{Nagaraju}] have proved common fixed point results based on this property. Later in 2013, Karapinar et al. \cite{Karapinar} gave the notion of this property in symmetric spaces for the two pairs of transformations. In this paper, we establish two common fixed point results for four transformations in vector $S$-metric spaces $(\Re, S, V)$ by using the notion of weakly compatibility (WC) and occasionally weakly compatibility (OWC).\\ \\
\textbf{\emph{Definition 1.1}} \cite{Ali}
Let $k, \rho:K \rightarrow K$ be two transformations. Then $k$ and $\rho$ are said to have a common fixed point $\alpha \in K$ if $k(\alpha)=\rho(\alpha)=\alpha$. \\ \\
\textbf{\emph{Definition 1.2}} \cite{Ali}
Let $k, \rho:K \rightarrow K$ be two transformations. Then $k$ and $\rho$ are said to have a coincidence point $\alpha \in K$ if $k(\alpha)=\rho(\alpha)=\beta$ and $\beta$ is called a point of coincidence of $k$ and $\rho$. \\ \\
\textbf{\emph{Definition 1.3}}\cite{Pooja}
A relation $\preceq$ on a set $\complement$  is a partial order if it satisfies the conditions mentioned below:
\begin{itemize}
\item[(a)]$ \flat_1 \preceq \flat_1  $ \hspace{7.4cm} (reflexive)
\item[(b)]$\flat_1\preceq \flat_2  $ and $ \flat_2\preceq \flat_1  $ implies $\flat_1=\flat_2$                                                                                                                 \hspace{2.3cm} (anti-symmetry)
\item[(c)]$ \flat_1\preceq \flat_2  $ and $ \flat_2\preceq \flat_3  $\,\,implies $ \flat_1\preceq \flat_3$                                                                           
\hspace{2.4cm} (transitivity)
\end{itemize}
$\forall$ $\flat_1, \flat_2, \flat_3 \in \complement$  . The pair $(\complement, \preceq)$ is known as partially ordered set.\\ \\
\textbf{\emph{Definition 1.4}}\cite{Aliprantis}
A poset is called a lattice if each set with two elements has an infimum and a supremum.\\ \\
\textbf{\emph{Definition 1.5}}\cite{Aliprantis}
An ordered linear space where the ordering is lattice is called vector lattice. This is also called Riesz space.\\ \\
\textbf{\emph{Definition 1.6}}
Let $V$ be a vector lattice and $\Re$ be a nonvoid set. Then vector $S$-metric is a transformation $S:\Re\times \Re\times \Re\rightarrow V$ on $\Re$ if it satisfies the conditions mentioned below:
\begin{itemize}
\item [(a)] $S(\flat_1, \flat_2, \flat_3) \succeq 0$,
\item[(b)]$S(\flat_1, \flat_2, \flat_3) = 0$ iff $\flat_1 = \flat_2 = \flat_3$,
\item[(c)]$S(\flat_1, \flat_2, \flat_3) \preceq S(\flat_1, \flat_2, \alpha)+S(\flat_2, \flat_2, \alpha)+S(\flat_3, \flat_3, \alpha)$,
\end{itemize}
\,\,for all  $\flat_1, \flat_2, \flat_3,\alpha \in \Re$.\\
The triplet $(\Re, S, V)$ is called vector $S$-metric space.\\ \\
\textbf{\emph{Lemma 1.7}}
If $(\Re, S, V)$ is a vector $S$-metric space, then  \[S(\vartheta, \vartheta, \mu)=S(\mu, \mu, \vartheta) \,\,\,\,\,\,\,\forall \mu, \vartheta \in \Re.\]
\textbf{\emph{Proof.}} Let  $\mu, \vartheta \in \Re$, then
\begin{eqnarray}
S(\vartheta, \vartheta, \mu) &\preceq &S(\vartheta, \vartheta, \vartheta)+S(\vartheta, \vartheta, \vartheta)+S(\mu, \mu, \vartheta)\\
&=&S(\mu, \mu, \vartheta) \nonumber\\
S(\mu, \mu, \vartheta)& \preceq& S(\mu, \mu, \mu)+S(\mu, \mu, \mu)+S(\vartheta, \vartheta, \mu)\\
&=&S(\vartheta, \vartheta, \mu) \nonumber
\end{eqnarray}
By (1) and (2), we get $S(\vartheta, \vartheta, \mu)=S(\mu, \mu, \vartheta)$.\\
\\ 
\textbf{\emph{Definition 1.8}} Let $(\Re, S, V)$ and $(\Re', S', V')$ be two vector $S$-metric spaces. A function $K:(\Re, S, V) \rightarrow (\Re', S', V')$ is continuous at $\alpha \in \Re$ if $\forall \,\,\,\,\langle \hslash_\flat\rangle \in \Re$ with $\hslash_\flat \xrightarrow{S,V} \alpha$ then $K(\hslash_\flat) \xrightarrow{S,V} K(\alpha)$.\\ \\ 
\textbf{\emph{Definition 1.9}} \cite{Shahzad}
Let $(\Re, S, V)$ be a vector $S$-metric space and  $P$ and $Q$ are two transformations in a vector $S$-metric space. Then the transformations are called weakly compatible (WC) if $PQ\alpha=QP\alpha$ whenever $P\alpha=Q\alpha$ for $\alpha \in \Re$. \\ \\
\textbf{\emph{Definition 1.10}}\cite{Shahzad}
Let $(\Re, S, V)$ be a vector $S$-metric space, $P$ and $Q$ are two transformations in a vector $S$-metric space. Then the transformations are called occasionally weakly compatible (OWC) if $PQ\alpha=QP\alpha$ for some $\alpha$ where $\alpha$ is coincidence point of $P$ and $Q$.\\ \\
\textbf{\emph{Remark 1.11}}
Every pair of (WC) transformations is (OWC) but converse not true.\\ \\
\textbf{\emph{Example1.12}} Let $\Re=[0, 5]$. Define the transformations $P, Q:\Re \rightarrow \Re$ by \[P(\xi)=\xi^2, \,\,\,Q(\xi)=2\xi \,\,\,\,\,\,\,\forall\xi\in \Re. \] Then 0 and 2 are coincidence points. Further $P(Q(0))=Q(P(0))$ and $P(Q(2))\neq Q(P(2))$. Then $(P, Q) $ is OWC but not WC.
\\ \\
\textbf{\emph{Definition 1.13}} \cite{Sintunavarat}
 Let $P, Q:\Re \rightarrow\Re$ be two transformations in a vector $S$-metric space $(\Re, S, V)$. Then the transformations are said to satisfy $CLR_Q$ (common limit in range of $Q$) property if $\exists \,\,\,\, \langle \hslash_n \rangle \in \Re$ such that 
 \[\lim_{n\rightarrow \infty }P\hslash_n= \lim_{n\rightarrow \infty }Q\hslash_n=Q\alpha \,\,\,where \,\,\, \alpha\in \Re\]\\
 \textbf{\emph{Example1.14}}
 Let $\Re\in [0, 10]$. Define the transformations $P$, $Q$ on $\Re$ as  \[P\hslash=\hslash+4,  \,\,\,\,\,Q\hslash=3\hslash \,\,\,\,\,\, \forall \hslash\in \Re.\] Define $\langle \hslash_n\rangle = \langle 2+\dfrac{1}{3n}\rangle$. Then \[\lim_{n\rightarrow \infty}P\hslash_n=\lim_{n\rightarrow \infty }Q\hslash_n=6=Q(2)\in \Re.\] Thus $P$ and $Q$ satisfy $CLR_Q$ property.\\ \\
\textbf{\emph{Definition 1.15}} \cite{Patel}
Let $(K, A)$ and $(P, Q)$ are two pairs of transformations and $(\Re, S, V)$ be vector $S$-metric space. Then $(K, A)$ and $(P, Q)$ are said to satisfy $CLR_{(A, Q)}$ property if $\exists \,\,\,\,\, \langle \hslash_n\rangle, \langle \vartheta_n\rangle \in \Re$ such that 
  \[\lim_{n\rightarrow \infty }K\hslash_n= \lim_{n\rightarrow \infty }A\hslash_n=\lim_{n\rightarrow \infty }P\vartheta_n= \lim_{n\rightarrow \infty }Q\vartheta_n=\alpha  \,\,\,\, with\,\,\, \alpha=Af=Q\ell\]where  $\alpha, f, \ell\in \Re$.\\ \\
\textbf{\emph{Example 1.16}}
 Let $\Re\in [0, 2]$. Define the transformations $K$, $A$, $P$ and $Q$ on $\Re$ as
 \[K\hslash=\dfrac{\hslash}{4}, \,\,\,\,\, A\hslash=\hslash-\dfrac{3}{28},\,\,\, \,\,\, P\hslash=\dfrac{\hslash}{7}, \,\,\,\,\,\, Q\hslash=\hslash-\dfrac{6}{28} \,\,\,\,\,\,\,\,\, \hslash\in \Re.\] Define $\langle \hslash_n\rangle, \langle \vartheta_n\rangle \in\Re$ such that $\langle \hslash_n\rangle=\langle \dfrac{1}{7} +\dfrac{1}{n} \rangle, \,\,\,\,\langle \vartheta_n\rangle=\langle\dfrac{1}{4}+\dfrac{1}{2n}\rangle\,\,\,\,n\in \mathbb{N}$. Then \[\lim_{n\rightarrow \infty} K\hslash_n=\lim_{n\rightarrow \infty}A\hslash_n= \lim_{n\rightarrow \infty}P\vartheta_n=\lim_{n\rightarrow \infty} Q\vartheta_n=\dfrac{1}{28} \] with $\dfrac{1}{28}=A(\dfrac{1}{7})=Q(\dfrac{1}{4})$ where    $\dfrac{1}{28}, \dfrac{1}{7}, \dfrac{1}{4}\in \Re$.\\ \\
 \\ 
 \textbf{\emph{Remark 1.17}}
If $K=P$ and $A=Q$ then CLR$_{(A, Q)}$ property reduces to $CLR_A$ property. \\ \\  
\section{Main Results}
\textbf{Theorem 2.1} Let $H, K, P$ and $Q$ are four transformations of a vector $S$-metric space $(\Re, S, V)$ satisfying the conditions mentioned below: 
\begin{itemize}
\item[(a)] $(H, P)$ and $(K, Q)$ satisfy CLR$_{(P, Q)}$ property
\item[(b)] $(H, P)$ and $(K, Q)$ are WC
\item[(c)]\begin{eqnarray}\label{eq1} 
S(H\hslash, H\vartheta, K\eta)\preceq \rho\, U(\hslash, \vartheta, \eta) \,\,\,\,\, \forall \hslash, \vartheta, \eta \in \Re \,\,\, and \,\,\,\rho\in [0, 1)
\end{eqnarray}
where 
\begin{eqnarray}
U(\hslash, \vartheta, \eta))&=& max \Big[S(P\hslash, P\vartheta, Q\eta), S(P\hslash, P\hslash, H\hslash ), S(Q\eta, Q\eta, K\eta),\nonumber\\&& \hspace{1cm}\dfrac{1}{2}\{S(P\hslash, P\vartheta, K\eta)+S(Q\eta, Q\eta, H\hslash)\}\Big]\nonumber
\end{eqnarray}
\end{itemize} 
then transformations $H, K, P$ and $Q$ have common fixed point in $\Re$ which is unique.\\
\textbf{Proof} Since $(H, P)$ and $(K, Q)$ satisfy CLR$_{(P, Q)}$ property, $\exists$  $\langle \hslash_n \rangle,\, \langle\vartheta_n\rangle\in \Re$ such that 
\[\lim_{n\rightarrow \infty }H\hslash_n= \lim_{n\rightarrow \infty }P\hslash_n= \lim_{n\rightarrow \infty }K\vartheta_n = \lim_{n\rightarrow \infty }Q\vartheta_n=\xi\]
with $\xi=Pa=Qb$ for some $\xi, a, b\in \Re$.\\
We assert that $Pa=Ha$.\\
By \eqref{eq1}, we have 
\begin{eqnarray}\label{eq2}
S(Ha, Ha, K\vartheta_n)\preceq \rho\, U(a, a, \vartheta_n)
\end{eqnarray}
where 
\begin{eqnarray}
U(a, a, \vartheta_n)&=& max \Big[S(Pa, Pa, Q\vartheta_n), S(Pa, Pa, Ha ), S(Q\vartheta_n, Q\vartheta_n, K\vartheta_n),\nonumber\\&& \hspace{1cm}\dfrac{1}{2}\{S(Pa, Pa, K\vartheta_n)+S(Q\vartheta_n, Q\vartheta_n, Ha)\}\Big]\nonumber\\
U(a, a, \vartheta_n)&\preceq& max \Big[S(Pa, Pa, Q\vartheta_n), S(Pa, Pa, Ha ), 2S(Q\vartheta_n, Q\vartheta_n, \xi)+\nonumber\\&& \hspace{1cm} S(K\vartheta_n, K\vartheta_n, \xi),\dfrac{1}{2}\{S(Pa, Pa, K\vartheta_n)+S(Q\vartheta_n, Q\vartheta_n, Ha)\}\Big]\nonumber\\
\lim_{n\rightarrow \infty }U(a, a, \vartheta_n)&\preceq& max \Big[\lim_{n\rightarrow \infty }S(Pa, Pa, Q\vartheta_n), \lim_{n\rightarrow \infty }S(Pa, Pa, Ha ), \nonumber\\&& \hspace{1cm}\lim_{n\rightarrow \infty }(2S(Q\vartheta_n, Q\vartheta_n, \xi)+ S(K\vartheta_n, K\vartheta_n, \xi)),\nonumber\\&& \hspace{1cm}\dfrac{1}{2}\{\lim_{n\rightarrow \infty }(S(Pa, Pa, K\vartheta_n)+S(Q\vartheta_n, Q\vartheta_n, Ha))\}\Big]\nonumber
\end{eqnarray}
\begin{eqnarray}
&\preceq& max \Big[S(\xi, \xi, \xi), S(\xi, \xi, Ha), 3S(\xi, \xi, \xi), \dfrac{1}{2}\{S(\xi, \xi, \xi)\nonumber\\&& \hspace{1cm}+S(\xi, \xi, Ha)\}\Big]\nonumber\\
&\preceq& max \Big[0, S(\xi, \xi, Ha),0, \dfrac{1}{2} S(\xi, \xi, Ha) \Big]\nonumber\\
&\preceq& S(\xi, \xi, Ha)\nonumber
\end{eqnarray}
By using $S(\vartheta, \vartheta, \mu)=S(\mu, \mu, \vartheta)$, we get
\[\lim_{n\rightarrow \infty }U(a, a, \vartheta_n)\preceq S(Ha, Ha, \xi). \hspace{7cm}\]
Taking limit $n \rightarrow \infty$ in \eqref{eq2}\\
\[S(Ha, Ha, \xi)\preceq \rho\lim_{n\rightarrow \infty }U(a, a, \vartheta_n)\preceq \rho\,S(Ha, Ha, \xi) \] 
which implies $S(Ha, Ha, \xi)=0$,  
i.e. $Ha=\xi.$\\
Thus $Ha=\xi=Pa.$\\
This shows that transformations $P$ and $H$ have a coincidence point $a$.\\
Since $(H, P)$ is WC, then we have $H(Pa)=P(Ha)$ and so $H\xi=P\xi.$\\
Now, we assert that $Kb=Qb.$ \\
By \eqref{eq1}, we have
\begin{eqnarray}\label{eq3}
S(H\hslash_n, H\hslash_n, Kb)\preceq \rho\, U(\hslash_n, \hslash_n, b)
\end{eqnarray}
where
\begin{eqnarray}
U(\hslash_n, \hslash_n, b)&=& max \Big[S(P\hslash_n, P\hslash_n, Qb), S(P\hslash_n, P\hslash_n, H\hslash_n ), S(Qb, Qb, Kb),\nonumber\\&& \hspace{1cm}\dfrac{1}{2}\{S(P\hslash_n, P\hslash_n, Kb)+S(Qb, Qb, H\hslash_n)\}\Big]\nonumber\\
&\preceq& max \Big[S(P\hslash_n, P\hslash_n, Qb), 2S(P\hslash_n, P\hslash_n,  \xi)+S(H\hslash_n, H\hslash_n, \xi),\nonumber\\&& \hspace{1cm} S(Qb, Qb, Kb),\dfrac{1}{2}\{S(P\hslash_n, P\hslash_n, Kb)+S(Qb, Qb, H\hslash_n)\}\Big]\nonumber
\\
\lim_{n\rightarrow \infty }U(\hslash_n, \hslash_n, b)&\preceq& max \Big[\lim_{n\rightarrow \infty }S(P\hslash_n, P\hslash_n, Qb), \lim_{n\rightarrow \infty }\{2S(P\hslash_n, P\hslash_n,  \xi)+\nonumber\\&& \hspace{1cm}S(H\hslash_n, H\hslash_n, \xi)\}, \lim_{n\rightarrow \infty }S(Qb, Qb, Kb),\nonumber\\&& \hspace{1cm}\dfrac{1}{2}\{\lim_{n\rightarrow \infty }(S(P\hslash_n, P\hslash_n, Kb)+S(Qb, Qb, H\hslash_n))\}\Big]\nonumber
\end{eqnarray}
\begin{eqnarray}
&\preceq& max \Big[S(\xi, \xi, \xi), 3S(\xi, \xi, \xi), S(\xi, \xi, Kb), \dfrac{1}{2}\{S(\xi, \xi, Kb)\nonumber\\&& \hspace{1cm}+S(\xi, \xi, \xi)\}\Big]\nonumber\\
&\preceq& max \Big[0, 0, S(\xi, \xi, Kb), \dfrac{1}{2} S(\xi, \xi, Kb) \Big]\nonumber\\
&\preceq& S(\xi, \xi, Kb).\nonumber
\end{eqnarray} 
Taking limit $n \rightarrow \infty$ in \eqref{eq3}\\
\[S(\xi, \xi, Kb)\preceq \rho\lim_{n\rightarrow \infty }U(\hslash_n, \hslash_n, b)\preceq \rho\, S(\xi, \xi, Kb) \] 
which implies $S(\xi, \xi, Kb)=0$, so
 $Kb=\xi.$\\
Thus $Kb=\xi=Qb.$\\
This shows that transformations $Q$ and $K$ have a coincidence point $b$ .\\
Since $(K, Q)$ is WC, then we have $K(Qb)=Q(Kb)$ and so $K\xi=Q\xi$.\\
We assert that $P\xi=H\xi=\xi$.\\
By \eqref{eq1}, we have
\begin{eqnarray}
S(H\xi, H\xi, Kb)\preceq \rho\,U(\xi, \xi, b)\nonumber
\end{eqnarray}

Since$  H\xi=P\xi $
\begin{eqnarray}\label{eq4}
S(P\xi, P\xi, Kb)\preceq \rho\,U(\xi, \xi, b)
\end{eqnarray}
where 
\begin{eqnarray}
U(\xi, \xi, b)&=& max \Big[S(P\xi, P\xi, Qb), S(P\xi, P\xi, H\xi), S(Qb, Qb, Kb),\nonumber\\&& \hspace{1cm}\dfrac{1}{2}\{S(P\xi, P\xi, Kb)+S(Qb, Qb, H\xi)\}\Big]\nonumber\\
U(\xi, \xi, b)&=& max\Big[S(P\xi, P\xi, \xi), S(P\xi, P\xi, P\xi), S(Kb, Kb, Kb),\nonumber\\&& \hspace{1cm}\dfrac{1}{2}\{S(P\xi, P\xi, \xi)+S(\xi, \xi, P\xi)\}\Big]\nonumber
\\
&=& max\Big[S(P\xi, P\xi, \xi), S(P\xi, P\xi, P\xi), S(Kb, Kb, Kb),\nonumber\\&& \hspace{1cm}\dfrac{1}{2}\{S(P\xi, P\xi, \xi)+S(P\xi, P\xi, \xi)\}\Big]\nonumber\\
&=& max \Big[S(P\xi, P\xi, \xi), 0, 0, S(P\xi, P\xi, \xi) \Big]\nonumber\\
&=& S(P\xi, P\xi, \xi).\nonumber
\end{eqnarray}
Hence \eqref{eq4} becomes\\
\[S(P\xi, P\xi, \xi)\preceq \rho\,S(P\xi, P\xi, \xi) \] 
which implies $S(P\xi, P\xi, \xi)=0$, so
 $P\xi=\xi.$\\
Thus $H\xi=P\xi=\xi.$\\
Finally, we assert that $K\xi=Q\xi=\xi.$ \\
By \eqref{eq1}, we have
\begin{eqnarray}\label{eq5}
S(Ha, Ha, K\xi))\preceq \rho\,U(a, a, \xi)
\end{eqnarray}
where
\begin{eqnarray}
U(a, a, \xi)&=& max \Big[S(Pa, Pa, Q\xi), S(Pa, Pa, Ha ), S(Q\xi, Q\xi, K\xi),\nonumber\\&& \hspace{1cm}\dfrac{1}{2}\{S(Pa, Pa, K\xi)+S(Q\xi, Q\xi, Ha)\}\Big]\nonumber
\end{eqnarray}
\begin{eqnarray}
U(a, a, \xi)&=& max \Big[S(\xi, \xi, K\xi), S(\xi, \xi, \xi), S(K\xi, K\xi, K\xi), \dfrac{1}{2}\{S(\xi, \xi, K\xi)\nonumber\\&& \hspace{1cm}+S(K\xi, K\xi, \xi)\}\Big]\nonumber\\
&=& max \Big[S(\xi, \xi, K\xi), 0, 0, \dfrac{1}{2}\{S(\xi, \xi, K\xi)+S(\xi, \xi, K\xi)\}\Big]\nonumber\\
&=& max \Big[0, S(\xi, \xi, K\xi) \Big]\nonumber\\
&=& S(\xi, \xi, K\xi).\nonumber
\end{eqnarray} 
Hence \eqref{eq5} becomes\\
\[S(\xi, \xi, K\xi)\preceq \rho S(\xi, \xi, K\xi) \] 
which implies $S(\xi, \xi, K\xi)=0$, so
 $K\xi=\xi.$\\
Thus $K\xi=Q\xi=\xi$.\\
Hence $H\xi=P\xi=K\xi=Q\xi=\xi.$\\
This shows that transformations $H, K, P$ and $Q$ have a common fixed point $\xi$  .\\
Let  transformations $H, K, P$ and $Q$ have another common fixed point $\alpha $. Then \[H\alpha=P\alpha=K\alpha=Q\alpha=\alpha.\]\\ 
By \eqref{eq1}, we have
\begin{eqnarray}\label{eq6}
S(H\xi, H\xi, K\alpha)\preceq \rho\,U(\xi, \xi, \alpha)
\end{eqnarray}
where
\begin{eqnarray}
U(\xi, \xi, \alpha)&=& max \Big[S(P\xi, P\xi, Q\alpha), S(P\xi, P\xi, H\alpha ), S(Q\alpha, Q\alpha, K\alpha),\nonumber\\&& \hspace{1cm}\dfrac{1}{2}\{S(P\xi, P\xi, K\alpha)+S(Q\alpha, Q\alpha, H\xi)\}\Big]\nonumber
\end{eqnarray}
\begin{eqnarray}
U(\xi, \xi, \alpha)&=& max \Big[S(\xi, \xi, \alpha), S(\xi, \xi, \alpha), S(\alpha, \alpha, \alpha), \dfrac{1}{2}\{S(\xi, \xi, \alpha)\nonumber\\&& \hspace{1cm}+S(\alpha, \alpha, \xi)\}\Big]\nonumber\\
&=& max \Big[S(\xi, \xi, \alpha), S(\xi, \xi, \alpha), 0, \dfrac{1}{2}\{S(\xi, \xi, \alpha)+S(\xi, \xi, \alpha)\}\Big]\nonumber\\
&=& max \Big[0, S(\xi, \xi, \alpha) \Big]\nonumber\\
&=& S(\xi, \xi, \alpha).\nonumber
\end{eqnarray} 
Hence \eqref{eq6} becomes\\
\[S(\xi, \xi, \alpha)\preceq \rho\,S(\xi, \xi, \alpha) \] 
which implies $S(\xi, \xi, \alpha)=0$, so
$\alpha=\xi.$\\
Hence transformations $H, K, P$ and $Q$ have common fixed point in $\Re$ which is unique.\\\\
%\textbf{Corollary 2.2} Let $H$ and $K$ are four self-maps of a vector $S$-metric space $(\Re, S, V)$ satisfy \eqref{eq1} and also if 
%\begin{itemize}
%\item[(1)] $H$ and $K$ satisfy CLR$_K$ property
%\item[(2)] $H$ and $K$ are WC
%\item[(3)]\begin{eqnarray}\label{eq1} 
%S(H\hslash, H\vartheta, K\eta)\preceq \rho\, U(\hslash, \vartheta, \eta) \,\,\,\,\, \forall \hslash, \vartheta, \eta \in \Re 
%\end{eqnarray}
%where 
%\begin{eqnarray}
%U(\hslash, \vartheta, \eta))&=& max \Big[S(H\hslash, H\vartheta, K\eta), \dfrac{1}{2}\{S(H\hslash, H\vartheta, K\eta)+S(K\eta, K\eta, H\hslash)\}\Big]\nonumber
%\end{eqnarray}
%\end{itemize} 
%then the mappings $H$ and $K$ have common fixed point in $\Re$ which is unique.\\
\textbf{Lemma 2.2}\cite{Harish} Let $\Re$ be a set, $P$ and $Q$ are OWC transformations of $\Re$. If $\omega=P\alpha=Q\alpha$ for some $\alpha\in \Re$ where $P$ and $Q$ have $\omega$ as point of coincidence which is unique then $\omega$ is common fixed point of $P$ and $Q$ which is unique.\\
\textbf{Proof:} The proof of this lemma can be shown as in the proof of \cite{Harish}.\\ \\
\textbf{Theorem 2.3} Let $H, K, P$ and $Q$ are four transformations of a vector $S$-metric space $(\Re, S, V)$ satisfying the condition mentioned below
\begin{eqnarray}\label{eq7}
S(H\hslash, H\vartheta, K\eta)\preceq \rho \,U(\hslash, \vartheta, \eta)\,\,\,\,\,\,\,\,\,\forall\,\,\, \hslash, \vartheta, \eta \in \Re\,\,\,and \,\,\,\rho\in [0, 1)
\end{eqnarray}
where 
%\begin{eqnarray}
\[U(\hslash, \vartheta, \eta)=max\{S(P\hslash, P\vartheta, Q\eta), S(P\hslash, P\hslash, H\hslash ), S(Q\eta, Q\eta, K\eta),S(Q\eta, Q\eta, H\hslash)\}\]
%\end{eqnarray}
and $(H, P)$ and $(K, Q)$ be OWC. Then transformations $H, K, P$ and $Q$ have common fixed point in $\Re$ which is unique.\\
\textbf{Proof} Since $(H, P)$ and $(K, Q)$ be OWC, then $\gamma, \delta\in \Re$ such that $H\gamma=P\gamma=\gamma_1$ and $H(P\gamma)=P(H\gamma)$ aand $K\delta=Q\delta=\delta_1$ and $K(Q\delta)=Q(K\delta).$\\
Firstly, we assert that $H\gamma=K\delta$.\\
By \eqref{eq7}, we have
\begin{eqnarray}\label{eq8}
S(H\gamma, H\gamma, K\delta)\preceq \rho\,U(\gamma, \gamma, \delta)
\end{eqnarray}
where 
\begin{eqnarray}
U(\gamma, \gamma, \delta)&=&max \{S(P\gamma, P\gamma, Q\delta), S(P\gamma, P\gamma, H\gamma), S(Q\delta, Q\delta, K\delta),\nonumber\\&& \hspace{1cm}S(Q\delta, Q\delta, H\gamma)\}\nonumber
\end{eqnarray}
\begin{eqnarray}
&=&max \{S(H\gamma, H\gamma, K\delta), S(H\gamma, H\gamma, H\gamma), S(K\delta, K\delta, K\delta),\nonumber\\&& \hspace{1cm} S(K\delta, K\delta, H\delta)\}\nonumber
\\
&=&max \{S(H\gamma, H\gamma, K\delta), S(H\gamma, H\gamma, H\gamma), S(K\delta, K\delta, K\delta),\nonumber\\&& \hspace{1cm}S(H\gamma, H\gamma, K\delta)\}\nonumber\\
&=&max \{S(H\gamma, H\gamma, K\delta),0, 0, S(H\gamma, H\gamma, K\delta)\}\nonumber\\
&=&\hspace{1cm} S(H\gamma, H\gamma, K\delta)\nonumber
\end{eqnarray}
Hence \eqref{eq8} becomes
\begin{eqnarray}
S(H\gamma, H\gamma, K\delta)\preceq \rho\,S(H\gamma, H\gamma, K\delta) \nonumber
\end{eqnarray}
which implies $S(H\gamma, H\gamma, K\delta)=0$, so $H\gamma=K\delta$.\\
Thus $H\gamma=P\gamma=K\delta=Q\delta$.\\
Let $\gamma_1\in\Re$ such that $H\gamma_1=P\gamma_1$.\\
By  \eqref{eq7}, we have
\begin{eqnarray}\label{eq9}
S(H\gamma_1, H\gamma_1, K\delta)\preceq \rho\,U(\gamma_1, \gamma_1, \delta)
\end{eqnarray}
where 
\begin{eqnarray}
U(\gamma_1, \gamma_1, \delta)&=&max \{S(P\gamma_1, P\gamma_1, Q\delta), S(P\gamma_1, P\gamma_1, H\gamma_1 ), S(Q\delta, Q\delta, K\delta),\nonumber\\&& \hspace{1cm}S(Q\delta, Q\delta, H\gamma_1)\}\nonumber\\
&=&max \{S(H\gamma_1, H\gamma_1, K\delta), S(H\gamma_1, H\gamma_1, H\gamma_1), S(K\delta, K\delta, K\delta),\nonumber\\&& \hspace{1cm} S(K\delta, K\delta, H\gamma_1)\}\nonumber\\
&=&max \{S(H\gamma_1, H\gamma_1, K\gamma), S(H\gamma_1, H\gamma_1, H\gamma_1), S(K\delta, K\delta, K\delta),\nonumber\\&& \hspace{1cm}S(H\gamma_1, H\gamma_1, K\delta)\}\nonumber\\
&=&max \{S(H\gamma_1, H\gamma_1, K\delta),0, 0, S(H\gamma_1, H\gamma_1, K\delta)\}\nonumber\\
&=&\hspace{1cm} S(H\gamma_1, H\gamma_1, K\delta)\nonumber
\end{eqnarray}
Hence \eqref{eq9} becomes
\begin{eqnarray}
S(H\gamma_1, H\gamma_1, K\delta)\preceq \rho S(H\gamma_1, H\gamma_1, K\delta)\nonumber
\end{eqnarray}
which implies $S(H\gamma_1, H\gamma_1, K\delta)=0$, so  $H\gamma_1=K\delta$.\\
Therefore $H\gamma_1=P\gamma_1=K\delta=Q\delta$.\\
Thus $H\gamma_1=H\gamma=P\gamma=P\gamma_1$.\\
Hence $\omega=H\gamma=P\gamma$ implies transformations $H$ and $P$ have unique point of coincidence.\\
By using lemma 2.2, transformations $H$ and $P$ have unique common fixed point, say $\xi$. Also in similiar way transformations  $K$ and $Q$ have unique common fixed point, say $\xi_1$.\\
Now we assert that $\xi=\xi_1$. \\
By \eqref{eq7}, we have
\begin{eqnarray}
S(H\xi, H\xi, K\xi_1)\preceq \rho\, U(\xi, \xi, \xi_1) \nonumber
\end{eqnarray}
By using $H\xi=P\xi=\xi$ and $K\xi_1=Q\xi_1=\xi_1$, we have
\begin{eqnarray}\label{eq10}
S(\xi, \xi, \xi_1)\preceq \rho\,U(\xi, \xi, \xi_1)
\end{eqnarray}
where 
\begin{eqnarray}
U(\xi, \xi, \xi_1)&=&max \{S(P\xi, P\xi, Q\xi_1), S(P\xi, P\xi, H\xi ), S(Q\xi_1, Q\xi_1, K\xi_1),\nonumber\\&& \hspace{1cm}S(Q\xi_1, Q\xi_1, H\xi)\}\nonumber\\
&=&max \{S(\xi, \xi, \xi_1), S(\xi, \xi, \xi), S(\xi_1, \xi_1, \xi_1), S(\xi_1, \xi_1, \xi)\}\nonumber\\
&=&max \{S(\xi, \xi, \xi_1), 0, 0, S(\xi_1, \xi_1, \xi)\}\nonumber\\
&=&max \{S(\xi, \xi, \xi_1),0, 0, S(\xi, \xi, \xi_1)\}\nonumber\\
&=&\hspace{1cm} S(\xi, \xi, \xi_1)\nonumber
\end{eqnarray}
Hence \eqref{eq10} becomes
\begin{eqnarray}
S(\xi, \xi, \xi_1)\preceq \rho S(\xi, \xi, \xi_1)\nonumber
\end{eqnarray}
which implies $S(\xi, \xi, \xi_1)=0$, so  $\xi=\xi_1$.\\
Thus transformations $H$, $K$, $P$ and $Q$ have common point in $\Re$ which is unique.\\ 
\end{document}